\documentclass[journal]{IEEEtran}
\usepackage{amsmath,amsfonts}
\usepackage{algorithmic}
\usepackage{algorithm}
\usepackage{array}
\usepackage[caption=false,font=normalsize,labelfont=sf,textfont=sf]{subfig}
\usepackage{textcomp}
\usepackage{stfloats}
\usepackage{url}
\usepackage{verbatim}
\usepackage{graphicx}
\usepackage{cite}
\usepackage{color} % for rev
\usepackage{hyperref}
\usepackage{scalerel}
\usepackage{tikz}
\usetikzlibrary{svg.path}

\definecolor{orcidlogocol}{HTML}{A6CE39}
\tikzset{
  orcidlogo/.pic={
    \fill[orcidlogocol] svg{M256,128c0,70.7-57.3,128-128,128C57.3,256,0,198.7,0,128C0,57.3,57.3,0,128,0C198.7,0,256,57.3,256,128z};
    \fill[white] svg{M86.3,186.2H70.9V79.1h15.4v48.4V186.2z}
                 svg{M108.9,79.1h41.6c39.6,0,57,28.3,57,53.6c0,27.5-21.5,53.6-56.8,53.6h-41.8V79.1z M124.3,172.4h24.5c34.9,0,42.9-26.5,42.9-39.7c0-21.5-13.7-39.7-43.7-39.7h-23.7V172.4z}
                 svg{M88.7,56.8c0,5.5-4.5,10.1-10.1,10.1c-5.6,0-10.1-4.6-10.1-10.1c0-5.6,4.5-10.1,10.1-10.1C84.2,46.7,88.7,51.3,88.7,56.8z};
  }
}

\newcommand\orcidicon[1]{\href{https://orcid.org/#1}{\mbox{\scalerel*{
\begin{tikzpicture}[yscale=-1,transform shape]
\pic{orcidlogo};
\end{tikzpicture}
}{|}}}}

\allowdisplaybreaks

\begin{document}

% \definecolor{r2}{rgb}{0,0,1} % blue
\definecolor{r2}{rgb}{0,0,0} % black
% \definecolor{r1}{rgb}{0,0,1} % blue
\definecolor{r1}{rgb}{0,0,0} % black

\title{Long-term Hydrothermal Bid-based Market Simulator}% with Case Studies in the Brazilian System}

\author{
Joaquim Dias Garcia$^{\orcidicon{0000-0002-7721-8564}}$,
Alexandre Street$^{\orcidicon{0000-0003-1569-030X}}$,
Mario Veiga Pereira %$^{\orcidicon{0000-0002-3771-1776}}$
% Joaquim Dias Garcia$^{\orcidicon{0000-0002-7721-8564}}$,~\IEEEmembership{Member,~IEEE,}
% Alexandre Street$^{\orcidicon{0000-0003-1569-030X}}$,~\IEEEmembership{Senior Member,~IEEE,}\\
% Mario Veiga Pereira,~\IEEEmembership{Fellow,~IEEE}
        % <-this % stops a space
\thanks{
Joaquim Dias Garcia (joaquim@psr-inc.com) is with PSR and LAMPS at PUC-Rio.
Alexandre Street (street@ele.puc-rio.br) is with LAMPS at PUC-Rio. Mario Veiga Pereira (mario@psr-inc.com) is with PSR.
The authors thank PSR's team for contributing with ideas, software, data, and infrastructure.
Joaquim Dias Garcia was partially supported by the Coordenação de Aperfeiçoamento de Pessoal de Nível Superior - Brasil (CAPES) - Finance Code 001.
Alexandre Street was also partially supported by FAPERJ and CNPq.
The work was partially funded by
the project P\&D ANEEL PD-00403-0050/2020, sponsored by ENGIE BRASIL ENERGIA S.A.
}% <-this % stops a space
%\thanks{Manuscript received April 19, 2021; revised August 16, 2021.}
}


% Remember, if you use this, you must call \IEEEpubidadjcol in the second

\maketitle

% number of chars:
\begin{abstract}
Simulating long-term hydrothermal bid-based markets considering strategic agents is a challenging task.
The representation of strategic agents considering intertemporal constraints within a stochastic framework brings additional complexity to the already difficult single-period bilevel, thus, non-convex, optimal bidding problem.
Thus, we propose a simulation methodology that effectively addresses these challenges for large-scale hydrothermal power systems.
We demonstrate the effectiveness of the framework through a case study with real data from the large-scale Brazilian power system.
In the case studies, we show the effects of market concentration in power systems and how contracts can be used to mitigate them. In particular, we show how market power might affect the current setting in Brazil.
The developed method can strongly benefit policymakers, market monitors, and market designers as simulations can be used to understand existing power systems and experiment with alternative designs.
\end{abstract}


\vspace{-6mm}
\section*{Nomenclature}
\label{sec:nomenclature}
% \printnomenclature

\vspace{-3mm}
\subsection*{Sets and Indices}

\begin{IEEEdescription}[\IEEEusemathlabelsep\IEEEsetlabelwidth{$MM$}]

\item[$I$]{Set of agents indexed by $i$.}
\item[$I_{-i}$]{Set of agents excluding agent $i$.}
\item[$I^M$]{Set of price-maker agents, indexed by $i$.}
\item[$I^T$]{Set of price-taker agents, indexed by $i$.}
\item[$J^G$]{Set of thermal plants indexed by $j$.}
\item[$J^H$]{Set of hydro plants indexed by $j$.}
\item[$J^L$]{Set of lags of autoregressive model, indexed by $l$.}
% \item[$J^M$]{Set of Markov states indexed by $j$.}
\item[$J^R$]{Set of renewable plants indexed by $j$.}
\item[$J^U(j)$]{Set of hydro plants upstream plant $j$, indexed by $y$.}
\item[$J^V$]{Set of vertices (points) representing a convex hull.}
\item[$S$]{Set of sampled scenario indices, indexed by $s$.}
\item[$T$]{Set of time indices, indexed by $t$.}
% \item[$K$]{SDDP Iteration index.}
% \item[{$ [K] $}]{Set of cut indices at iteration $K$.}
\item[{$J_i$}]{Subset of a set $J$ with the indices of elements that belong to agent $i$.}

\end{IEEEdescription}

Subscripts related to stage, $t$, and scenario, $s$, will be omitted for simplicity when they are not essential for the reader's understanding of sampling and chronological relations.

\vspace{-3mm}
\subsection*{Constants}

\begin{IEEEdescription}[\IEEEusemathlabelsep\IEEEsetlabelwidth{$MM$}]

\item[$\mathcal{A}$]{Vector of inflows of all hydros for all stages and sampled scenarios.}
\item[$C_j$]{Operating cost of thermal $j$.}
\item[${E}^{Q}_j$]{Energy quantity for an element $j$ of a convex hull.}
\item[${E}^{R}_j$]{Revenue for an element $j$ of a convex hull.}
\item[$G_j$]{Maximum generation of thermal $j$.}
\item[$P_i$]{Price offer of agent $i$.}
\item[$\mathcal{P}$]{Vector price offer of all agents for all stages and sampled scenarios.}
\item[$P^F$]{Forward contract price.}
\item[$Q_i$]{Quantity offer of agent $i$.}
\item[$\mathcal{Q}$]{Vector of quantity offer of all agents for all stages and sampled scenarios.}
\item[$Q^F$]{Forward contract quantity.}
\item[$\tilde{R}_j(\omega)$]{Maximum generation of renewable $k$ at scenario $\omega$.}
\item[$\mathcal{R}$]{Vector of maximum generation of all renewables for all stages and sampled scenarios.}
\item[$U_j$]{Maximum flow through turbine of hydro $j$.}
\item[$V_j$]{Maximum storage of hydro $j$.}

% \item[$\beta^k$]{Constant term of Benders cut $k$.}
% \item[$\gamma^k_j$]{Coefficient of state $v_j^{t+1}$ of Benders cut $k$.}
% \item[$\delta^k_{j,l}$]{Coefficient of state $a_j^{t+1-l}$ of Benders cut $k$.}
\item[$\tilde{\varepsilon}_i^t(\omega)$]{Inflow noise coefficient of hydro $j$, stage $t$ and scenario $\omega$.}
\item[$\phi_{j,l}$]{Inflow autoregressive coefficient of hydro $j$, lag $l$.}

\item[$\pi$]{Spot price.}
\item[$\mathit{\Pi}$]{Vector of spot prices for all stages and sampled scenarios.}

% \item[$M_{\mu|m}$]{Transition probability from state $\mu$ to state $m$.}
\item[$\mathcal{M}$]{Vector of transition probabilities matrices, for all stages.}

\end{IEEEdescription}

Indexing vectors in calligraphic ($\mathcal{P}, \mathcal{Q}$) by $i$ stands for the vector where elements belong to agent $i$ and by $-i$ stands for the vector where elements belong to all agents except $i$.
\vspace{-3mm}
\subsection*{Optimization Variables}

\begin{IEEEdescription}[\IEEEusemathlabelsep\IEEEsetlabelwidth{$MM$}]

\item[$a_j^t$]{Inflow at hydro $j$, stage $t$. $a^{[t]}$ stands for the vector of all inflows before stage $t$.}
\item[$\textbf{a}^{[t]}$]{Vector of inflows at all hydros, for all stages before stage $t$.}
\item[$e$]{Energy offer.}
\item[$g_j$]{Generation of thermal $j$.}
\item[$q$]{Bid quantity accepted in a dispatch.}
\item[$r_j$]{Generation of renewable $j$.}
\item[$u_j$]{Turbine flow at hydro $j$.}
\item[$v_j^{t}$]{Storage at hydro $j$, at the beginning of stage $t$, and at the end of stage $t-1$.}
\item[$\textbf{v}^{t}$]{Vector of storage at all hydros, at the beginning of stage $t$, and the end of stage $t-1$.}
\item[$z_j$]{Spill flow at hydro $j$.}
% \item[$\alpha$]{Epigraph variable for Benders cuts.}
\item[$\lambda_j$]{Convex hull value for vertex $j$.}

\end{IEEEdescription}

Indexing vectors in bold ($\textbf{a}, \textbf{v}$) by $i$ stands for the sub-vector where elements belong to agent $i$.

\vspace{-3mm}
\subsection*{Functions and Functionals}

\begin{IEEEdescription}[\IEEEusemathlabelsep\IEEEsetlabelwidth{$MM$}]

\item[{$\mathbb{E}_\omega[\cdot]$}]{Expected value over the random variable $\omega$.}
\item[$\tilde{B}(\cdot, \cdot)$]{Future cost as a function of states and uncertainty.}
\item[$\tilde{\Lambda}(\cdot, \cdot)$]{Revenue as a function of energy, $e$ and scenario $\omega$.}
\item[$\rho(\cdot)$]{Hydro generation as a function of turbine flow. It can also be a function of volume.}
\item[$|\cdot|$]{Cardinality of a set.}

\end{IEEEdescription}

\vspace{-3mm}
\section{Introduction}
\IEEEPARstart{H}{ydro} power is one of the most widely used energy sources around the globe, the most used renewable energy source responsible for over four thousand TWh per year according to the International Energy Agency \cite{iea2022}. Many countries rely on hydro plants for a meaningful share of their generation, but some countries have hydropower as their main energy source, for instance, Brazil, New Zealand, Colombia, and Norway.

Many countries were initially centrally operated in a cost-based market design. In this market design, the water values and the dispatch of all units are centrally calculated based on audited costs by the system operator. This cost-based market is still the case in the Brazilian, Chilean, Mexican, Vietnamese, and South Korean power markets, to mention a few.

The changes in the regulation of power systems in the last decades have led to liberalized bid-based power markets, including the hydrothermal ones \cite{ieee2019conversations}. Norway, Colombia, New Zealand, and many other countries have implemented this market design change. In the Brazilian case, since the 2000 crisis, many waves of proposals for changing from an audited cost to a bid-based market design have been made starting with the power sector revitalization committee \cite{mme2002revit}. More recently, in 2017, a public call for contributions named CP33 was made in this country with the leadership of the system planner (EPE) \cite{epe2017cp33}, and guidelines for market reforms were established. The study of a bid-based market design in Brazil constitutes one of the key aspects of the new guidelines. This led to a working group in 2019 \cite{ccee2019gt} conducted by the market operator (CCEE), with the participation of the power system operator (ONS), the system planner (EPE) and the energy ministry (MME). Finally, there is an ongoing discussion around a law project from 2021 \cite{camara2021law} that includes a market liberalization pathway.

The complexity in assessing the potential of market power abuse within the large-scale setting of the Brazilian system with realistic data has always constituted a barrier to further discussion (see \cite{ribeiro2023technical} for more details). The main difficulty is the simulation of hydrothermal bid-based markets that combines the challenge of assessing the opportunity costs of water and the interactions among players in a uniform-price-based system. 
The bidding problem for hydrothermal power systems was reviewed in \cite{steeger2014optimal}, which contrasts techniques required in two main dimensions: 1) the need to represent hydro plants as opposed to the purely thermal case; 2) the representation (or not) of strategic agents that can manipulate markets.

The operation of systems with large penetration of hydropower is very challenging due to the need to model long time horizons and random variables. The first challenge is due to hydro reservoirs, which are long-duration energy storage devices that create a strong time coupling in the problem because water can be used for cheap generation or saved for future use. Saving too much water might lead to spillages, while excessive usage might lead to energy deficits. The second challenge is due to the uncertainty of inflows and renewable generation (such as solar and wind). These two aspects give rise to an opportunity cost of using the water today, known as the water value.
Multi-stage stochastic Optimization (MSO) is the most used framework to handle this type of problem, and the key algorithm is Stochastic Dual Dynamic Programming, SDDP \cite{pereira1991multi}, initially motivated by the optimal operation of the centralized hydrothermal systems.

On the other hand, the interaction between agents in a market requires simulation tools from Game Theory \cite{fudenberg1991game}. Modeling a single agent capable of market power involves the usage of Stackelberg Games modeled through a Bilevel Optimization \cite{pozo2017basic} that is a powerful framework but leads to non-convex optimization problems. The interaction of multiple agents leads to even more complex models that usually rely on Nash Equilibrium on top of the previous Stackelberg Games \cite{barroso2006nash, gabriel2012complementarity}.

Due to all this complexity, previous research has focused on simplifying parts of the problem to simulate systems.
Multi-agent hydrothermal markets were simulated through equilibrium modeling in the work  \cite{bushnell2003mixed} that applied Equilibrium Problems with Equilibrium Constraints (EPEC) but only considered 7 stages. EPECs were also used in  \cite{almeida2013medium} that modeled one year of operation by was limited to 4 scenarios, while \cite{yanjie2020power} considers two stages and ten scenarios, \cite{moiseeva2018bayesian} considers four stages and 20 scenarios, \cite{cruz2016strategic} considers 24 stages and one scenario.
The early works \cite{scott1996modelling} and \cite{kelman2001market} apply variations of Dynamic Programming (DP) to consider more stages and scenarios but require simple models of each stage that have analytical solutions. \cite{villar2003hydrothermal} and  \cite{steeger2015strategic} also employ DP, but in a deterministic equilibrium case study with few stages, while \cite{loschenbrand2018hydro} used a modified DP in a case study with 15 scenarios.
SDDP was first used in this setting to model a single strategic agent where the remainder of the market is purely thermal was developed in the seminal work \cite{flach2010long} and subsequently by \cite{steeger2017dynamic,habibian2020multistage}.

As per the previously reported works, simulating multiple agents in a realistic hydrothermal bid-based market is exceptionally challenging. To that end, one must handle multiple decision-making stages under uncertainty and the interplay of various agents, leading to non-convexities. Consequently, most works strongly limit the number of stages, agents, or scenarios. Although some of the above methods can handle all these features, in theory, the simulations are done in small problems, preventing us from deriving relevant practical studies of real large hydrothermal power systems, such as the Brazilian one.

Therefore, the main contributions of this work are:

\begin{itemize}
    \item Developing a methodology that can handle, at the same time, multiple stages, agents, and uncertainty scenarios to enable realistic studies on large hydrothermal power systems. Such methodology will be based on solving multi-stage stochastic strategic bidding problems for each strategic agent, while the coupling between agents is achieved by an iterative procedure based on diagonalization. Although SDDP and Game Theory were previously applied in this area, the uncertainty modeling of unknown bids from competing agents and an effective simulator based on the diagonalization method are novel.
    \item Simulating and analyzing the competition of multiple price-maker agents in the large-scale Brazilian power system with real data. In this sense, we provide new numerical studies and insights that may help in current discussions on possible changes toward a bid-based market in Brazil. Within this context, we also derive quantitative results on how agents' market share and forward involvement can mitigate market power abuse. 
    
\end{itemize}

The remainder of this work is organized as follows: Section \ref{sec:fw} presents important fundamentals that will be necessary for the following sections. Section \ref{sec:ht} introduces basic notation describing a centralized model for hydrothermal power systems. Section \ref{sec:sa} describes the optimization of a single price-maker agent and highlights how to incorporate contracts. Section \ref{sec:ma} details an algorithm to combine the above models to simulate long-term hydrothermal bid-based markets in the presence of multiple strategic agents. Section \ref{sec:case} presents case studies to test the proposed algorithm. Section \ref{sec:br} presents the simulation of a real case of the Brazilian power system. Finally, Section \ref{sec:con} exposes the main conclusions.

\vspace{-3mm}
\section{Methodological background} \label{sec:fw}

\subsection{Electricity markets}\label{sec:em}

Electricity markets are usually divided into two main groups: centralized markets (also known as cost-based) and liberalized markets (also known as bid-based). In the following, we present two figures to help contrast these two designs. Other market designs are possible \cite{ribeiro2023technical}, and some systems might mix the two (to some extent), like Brazil, which is cost-based but allows for flexible demand bids (small compared to the total load).

\subsubsection{Centralized cost-based markets}\label{sec:cd} 
\begin{figure}[!t]
% \centering
\includegraphics[width=3.1in, trim={-1.0cm 0.35cm 0.cm 0.02cm}, clip]{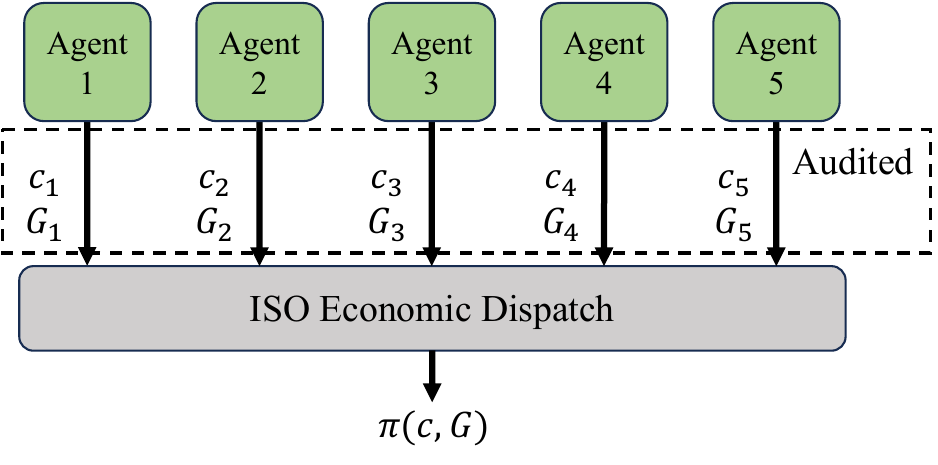}
\caption{Centralized cost-based markets.}
\label{fig:meth:cd}
\end{figure}

Their structure and information scheme are depicted in Figure \ref{fig:meth:cd}. In this design, agents have their costs, $c$, available capacity, $G$, and other relevant data audited by the system operator and/or regulator. With this information, the Independent System Operator (ISO) will optimize the system economic dispatch, which will lead to prices, $\pi$, that are, consequently, functions of $c$ and $G$. The main drawbacks are the burden of the auditing process and the lack of incentives for agents to improve and maintain their equipment and even to disclose information to the system operator. Centralization also creates the issue of single models and views, which, if mistaken or imprecise, can be questioned by agents and lead to a massive judicialization in the market. Thus, the decentralization of responsibility also supports the decentralization process towards a bid-based market design. 

\subsubsection{Liberalized bid-based markets}\label{sec:bb} 
\begin{figure}[!t]
% \centering
\includegraphics[width=3.0in, trim={0.cm 0.29cm 0.cm 0.05cm}, clip]{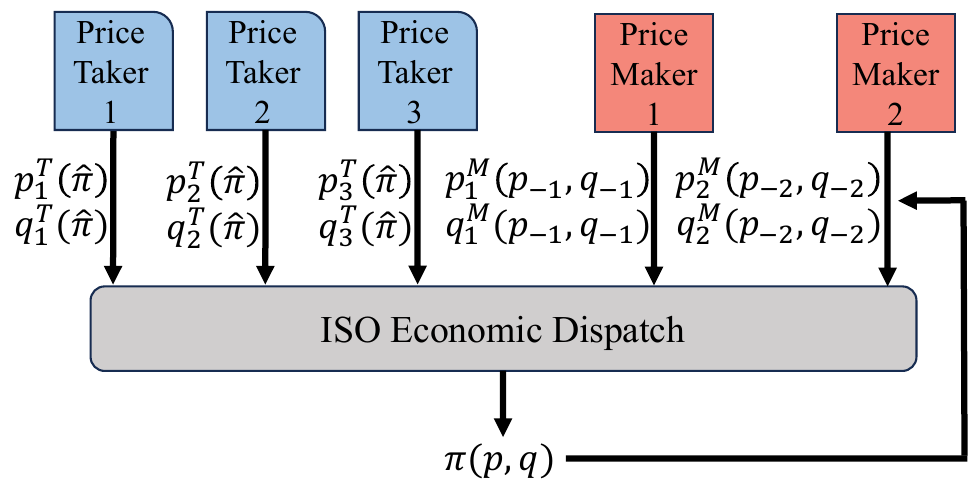}
\caption{Bid-based market.}
\label{fig:meth:bb}
\end{figure}

Their structure and information scheme is depicted in Figure \ref{fig:meth:bb}. In this case, there might be two main types of agents: price takers and price makers. Both send price, $p$, and quantity, $q$, bids to the ISO, which uses the bids to optimize the economic dispatch and define the spot price. Thus, in this market design, dispatch decisions and spot prices are functions of the offered pairs $p,q$. {\color{r2}Differences between these agent types will be described next.} %In the next paragraphs, we differentiate the two and explain the feedback arrow in the figure.

{\color{r2}
\subsection{Agent behavior in liberalized bid-based markets}\label{sec:ab}

%We describe the main differences between
The two types of agents introduced in Section \ref{sec:bb}, price takers and price makers (see Figure \ref{fig:meth:bb}), aim to optimize their bids and maximize profits, but using different approaches and methods.
}

\subsubsection{Price-taker agents}\label{sec:pt} These are non-strategic agents that are either too small to affect market prices by changing their operation or, despite being large, they make their offer assuming the price as given. In this context, we assume the latter for generality. Thus, this type of agent develops his bidding strategy, $(p,q)$, to maximize profits assuming the prices,  $\hat{\pi}$, as exogenous scenarios. Under mild conditions, markets with only price-taker thermal agents are solved in the literature with offers based on generators' variable costs and their total capacity \cite{gross2000generation,barroso2012overview}. Meanwhile, hydro agents should bid based on their water values, which can be calculated from the operator's perspective \cite{pereira1991multi} or from the agent's perspective \cite{gjelsvik1999algorithm}.

\subsubsection{Price-maker agents}\label{sec:pm} Are also named strategic agents and can be defined as the generators that can potentially change prices by changing their offers and, consequently, their operation. These agents are those that can potentially exercise market power abuse (the main drawback of a bid-based market). In this case, if a strategic agent considers the price response, $\pi(p,q)$, when developing its offers to maximize profits, and the offer actually changes the price, we say that the agent has exercised market power abuse. Hence, the feedback loop in Figure \ref{fig:meth:bb} highlights this possible relation between price marker bids that depend on their belief of the bid of other agents, $(p_{-i}, q_{-i})$, and the final prices.
In the purely thermal case,  \cite{barroso2006nash, fanzeres2020column} model the Nash equilibrium using Karush-Kuhn-Tucker optimality conditions (KKT) conditions of the market operation leading to a large-scale mixed-integer problem (MIP).

\vspace{-4.3mm}
\subsection{Equilibrium simulation}\label{sec:eq}

We present two equilibrium models and associated solution methods that will be instrumental for our developments.

\subsubsection{Stackelberg equilibrium}\label{sec:st} It is applicable when there is a single price-maker agent. In this case, the price maker is known as the leader agent, and the rest of the agents are the followers. The leader will attempt to optimize its revenue by strategically choosing a bid that will affect prices. This equilibrium is modeled with Bilevel Optimization \cite{pozo2017basic}.
It was extended to the hydro case by \cite{flach2010long} that embeds a bilevel program in the SDDP algorithm but convexifies each sub-problem to follow the convexity assumption of SDDP (see next Section \ref{sec:sddp}). A variant of \cite{flach2010long} was proposed in \cite{steeger2017dynamic} that handles the non-convexity with Lagrangian relaxation combined with SDDP. \cite{habibian2020multistage} uses a purely Lagrangian decomposition scheme to solve a multi-stage MIP for a price-maker demand side agent.

\subsubsection{Nash Equilibrium}\label{sec:nash} It applies to a more general setting with multiple price-maker agents. This state of equilibrium is defined as a set of offers for all agents such that no agent can improve the revenue by unilaterally deviating from this point \cite{fudenberg1991game, barroso2006nash}. Mathematically, it can be modeled as an Equilibrium Problem with Equilibrium Constraints (EPEC) derived from the KKT conditions of each agent revenue maximization problem. Analytical formulas can be found in special cases, as in \cite{scott1996modelling} and \cite{kelman2001market}, facilitating hydro players' multi-stage analysis through DP decomposition schemes. Another possibility is to apply MIP-based reformulations to single-stage equilibrium problems in the lines of \cite{barroso2006nash}:  \cite{steeger2015strategic} and \cite{loschenbrand2018hydro} model the hydro case by solving MIP within DP.
For small hydro problems with few scenarios and stages, \cite{bushnell2003mixed,almeida2013medium}, EPECs (and multi-stage EPECS) can be solved by specialized solvers such as PATH \cite{dirkse1995path}. 

\subsubsection{The diagonalization method}\label{sec:diag} Is sometimes referred to as nonlinear Gauss-Seidel or agent-based simulation, and it is especially well documented in \cite{gabriel2012complementarity} (see also \cite{devine2023strategic} and \cite{hori2019gauss}). In general terms, the method aims to solve an equilibrium problem by optimizing one agent at a time. Hence, the first agent optimizes its strategy with the strategy of the other agents fixed. Then, the second agent is optimized in the same way, but now, considering the strategy of the first agent fixed to the previously obtained value. If the first agent's strategy has changed, the second's current strategy may also change. An iteration of the algorithm is completed once all agents have been optimized. The iterations continue until the strategies from all agents do not change anymore. This method has been applied to very limited-sized hydrothermal equilibrium problems in
\cite{almeida2013medium},
\cite{yanjie2020power},
\cite{cruz2016strategic},
\cite{villar2003hydrothermal}, as mentioned in the introduction.

\vspace{-4.3mm}
\subsection{Multi-stage stochastic programming and SDDP}\label{sec:sddp}

MSO is a framework to handle problems that are coupled in time and depend on random variables. We can formulate the linear programming version through its Bellman recursion, which depicts the problem being solved in each stage and how they are coupled through state variables $\textbf{x}^{t}$:
\begin{align}
   &\tilde{B}^{t}\big(\textbf{x}^{t}, \omega^{t}\big) = \notag \\
   &\min_{\textbf{y}^{t}\in Y(\omega^{t}), \textbf{x}^{t+1}} c(\omega^{t})^{\top} \textbf{y}^{t} + \mathbb{E}_{\omega^{t+1}}\left[\tilde{B}^{t+1}\big(\textbf{x}^{t+1}, \omega^{t+1}\big) \right]\label{mod:mso:beg}\\
   &s.t.   \  \ \textbf{x}^{t+1} + A_x(\omega^t)\textbf{x}^{t} + A_y(\omega^t)\textbf{y}^{t} = b(\omega^t)\label{mod:mso:end}
\end{align}
SDDP is a fundamental algorithm for effectively solving MSO and is applicable as long as they follow two main assumptions: 1) convexity of stage problems (which is trivial in the above linear case); 2) stagewise independence, that is, for all $t$, $\omega^t$ is independent of $\omega^\tau$ if $\tau<t$.

{\color{r2}To apply the SDDP to the above problem, we need to reformulate it. The reformulation consists of replacing the future cost function (FCF) -- the second term in the objective function \eqref{mod:sddp:beg} -- by its piece-wise linear formulation, which is represented by the epigraph variable, $\alpha$, and cuts \eqref{mod:sddp:end} as follows:}
%To apply SDDP to the above problem, we need to reformulate it. The reformulation consists of replacing the future cost function (FCF), the second term in the objective function \eqref{mod:sddp:beg}, by a piece-wise linear form represented by the epigraph variable, $\alpha$, and the cuts \eqref{mod:sddp:end}:
% To apply SDDP to the above problem, we need to reformulate it so that the future cost function (FCF), the second term in the objective function \eqref{mod:sddp:beg}, in its piece-wise linear formulation:
%
\begin{align}
   \hspace{-1mm}\tilde{B}^{t}\big(\textbf{x}^{t}, \omega^{t}\big) =
   &\min_{\textbf{y}^{t}\in Y(\omega^{t}), \textbf{x}^{t+1}, \alpha} c(\omega^{t})^{\top} \textbf{y}^{t} + \alpha\label{mod:sddp:beg}\\
   &s.t.   \  \ \textbf{x}^{t+1} + A_x(\omega^t)\textbf{x}^{t} + A_y(\omega^t)\textbf{y}^{t} = b(\omega^t)\label{mod:sddp:bal}\\
   & \alpha \geq \beta^k+{\chi^k}^\top \textbf{x}^{t+1}, \quad \forall k \in \mathcal{K}\label{mod:sddp:end}
\end{align}
The above model considers a theoretical representation of the true model, as all the cuts %(represented by coefficients $\beta^k$ and $\chi^k$)
representing the FCF in \eqref{mod:sddp:end} are not known in advance. Therefore, the SDDP relaxes the set of all cuts, $\mathcal{K}$, to a subset $[K]$ that is iteratively updated. In this case, $\tilde{B}^{t}$ is also replaced by $\tilde{B}^{t}_K$ as it will only be a lower approximation of $\tilde{B}^{t}$ (at iteration $K$).

In a very high-level description, the SDDP algorithm starts sampling a scenario, $\omega^t$, for each $t\in T$, then solves the stage subproblems in chronological order, generating a feasible solution, this phase is called the forward step. After that, in the so-called backward step, problems are solved in the reverse order of time, generating Benders cuts to improve the representation of the FCF and, thus, propagating information from the future to the present. By averaging multiple feasible solutions for forward passes, we compute a statistical upper bound, while the solution of the first stage problem gives a deterministic lower bound.
These steps are repeated until a specified stopping criterion is reached. This process is also known as \textit{policy optimization} (or training). For algorithm details, the reader is directed to \cite{pereira1991multi, shapiro2011analysis}. After convergence is achieved, it is usual to proceed with a final forward pass in which we sample scenarios $s \in S$, and for each of them, we solve all stage subproblems in chronological order. This last procedure, known as \textit{simulation}, results in a solution for each primal and dual variable. In other words, we finish with values for each stage in $T$ and sample scenario in $S$ for all optimization variables.

In this section, we focused our presentation on the expected value case of MSO for simplicity. However, it is relevant to mention that the methodological developments considered in this paper hold for other risk measures, such as the conditional value at risk, that can be consistently represented in the SDDP scheme (see \cite{shapiro2013risk, rudloff2014time,dowson2022incorporating}).

\subsubsection{Uncertainty representation in SDDP with auto-regressive processes}\label{sec:ts-sddp} It is a technique to handle stagewise-dependent uncertainty. This technique requires that the uncertainty only affects the right-hand side (RHS) of \eqref{mod:sddp:bal}, $d(\omega)$, and it can be represented by an auto-regressive stochastic process with stagewise independent noise. In this case, we expand the state variables to include lagged observations of the uncertainty and add a constraint representing the evolution of the process in the model. Therefore, the final model will have the same form as \eqref{mod:sddp:beg}--\eqref{mod:sddp:end} but with more variables and constraints. The uncertainty will be restricted to the RHS of \eqref{mod:sddp:bal} and stagewise independent. This is the case of inflows modeled this way since the first appearance of SDDP \cite{pereira1991multi}.

\subsubsection{Uncertainty representation in SDDP with Markov processes}\label{sec:m-sddp} It is a second technique to handle stagewise-dependent uncertainty.
While the previous method effectively handles RHS uncertainty, it does not apply to the cases where stagewise dependent uncertainty affects objective or left-hand side coefficients (coefficients that multiply variables). To handle these cases, an alternative version of SDDP that combined ideas from SDDP and Stochastic Dynamic Programming (SDP) was developed (in the context of price-taker optimal hydro bidding) by \cite{gjelsvik1999algorithm}. More recently, this approach has been termed Markovian SDDP \cite{lohndorf2019modeling} as the time dependency of the uncertainty is handled by a Markov Chain model. For a given stage, $t$, and Markov state, $\mu$, the uncertainty is now considered independent of the past, but the dependency structure is preserved in the transition probabilities.
The reformulated model contains multiple FCFs, one FCF for each possible future Markov state, $m$, represented by the epigraph variable $\alpha_m$, and the respective cuts in \eqref{mod:sddpmk:end}:
\begin{align}
   &\tilde{B}^{t}_{\mu, K}\big(\textbf{x}^{t}, \omega^{t}_\mu\big) = \notag \\
   &\min_{\textbf{y}^{t}\in Y(\omega^{t}), \textbf{x}^{t+1}, \alpha_m} c(\omega^{t})^{\top} \textbf{y}^{t} + \sum_{m \in J^M} M_{\mu|m} \alpha_m\label{mod:sddpmk:beg}\\
   &s.t.   \  \ \textbf{x}^{t+1} + A_x(\omega^t)\textbf{x}^{t} + A_y(\omega^t)\textbf{y}^{t} = b(\omega^t)\label{mod:sddpmk:bal}\\
   & \alpha_m \geq \beta^k_m+{\chi^k_m}^\top \textbf{x}^{t+1}, \quad  \forall m \in J^M,k \in \mathcal{K}\label{mod:sddpmk:end}
\end{align}
Where $M_{\mu|m}$ is the transition probability from Markov state $m$ to $\mu$ and $J^M$ is the set of Markov states.

The uncertainty representation methods described above are not mutually exclusive. A model might represent some uncertainties with the auto-regressive reformulation and some with the Markovian.

\vspace{-2mm}
\section{Centralized long-term hydrothermal dispatch and SDDP} \label{sec:ht}

The hydrothermal power system dispatch is a very complex problem because it is a multi-stage stochastic optimization problem that includes many physical and policy-driven constraints. From now on, we will also refer to the market design where generation dispatch orders and spot prices are derived based on centralized long-term dispatch calculations as a centralized audited-cost-based market.

We present a simple yet general optimization model with all the main core features required in the market simulator proposed in this work.
We describe the problem as a Bellman recursion as in Section \ref{sec:sddp}. Therefore, \eqref{mod:ht:init}--\eqref{mod:ht:end} presents the objective function and constraints of a given stage, $t$, and random event $\omega^t$. Index $t$ is omitted from most variables when easily understood from context aiming at a lighter notation.

The first formula, \eqref{mod:ht:init}, states that the future cost of stage $t-1$, namely $\tilde{B}^{t}$, given the states $\{\textbf{v}^{t}, \textbf{a}^{[t-1]}\}$ (vector of water stored and observed past inflows in each reservoir) and the random event $\omega^t$, is defined as the minimization of the problem \eqref{mod:ht:init}--\eqref{mod:ht:end}. In this model, the objective function can be split into two pieces. The first piece is the immediate cost in the form of a thermal cost (load shedding accounted for as an additional artificial expensive generator), while the second piece is the expected value of the future cost $\tilde{B}^{t+1}$, where $\tilde{B}^{|T|+1}(\cdot) = 0$.
The equation \eqref{mod:ht:load} represents the load balance of the system. We have generation and energy flows on the left-hand side and demand on the right-hand side. Demand is considered deterministic to simplify the developments, but everything that follows can be trivially extended to the stochastic demand case. 
\eqref{mod:ht:mass} describes the water mass balance: storage at the end of stage $t$ equals the storage at the beginning plus the net sum of incoming water and outflows.
\eqref{mod:ht:hyd} constrains hydro storage, turbine flow, and spillage, while \eqref{mod:ht:ter} limits the thermal generation and \eqref{mod:ht:ren} limits the renewable generation that has a stochastic upper bound depending on events like sun and wind.
Finally, \eqref{mod:ht:end} describes the inflow autoregressive stochastic process as in Section \ref{sec:ts-sddp} and \cite{pereira1991multi,infanger1996cut}.
\begin{align}
   & \tilde{B}^{t}\big(\textbf{v}^{t}, \textbf{a}^{[t-1]}, \omega^{t}\big) = \notag\\
   &\min \ \ \sum_{j \in J^G} C_j g_j + \mathbb{E}_{\omega^{t+1}}\left[\tilde{B}^{t+1}\big(\textbf{v}^{t+1}, \textbf{a}^{[t]}, \omega^{t+1}\big) \right]\label{mod:ht:init}\\
   &s.t.   \notag  \\
   & \sum_{j\in J^G} g_j + \sum_{j\in J^H} \rho_j(u_j) + \sum_{j\in J^R} r_j 
   = {D}  \label{mod:ht:load} \\
   & v^{t+1}_j = v^t_j - u_j - z_j + \sum_{n\in J^U(j)} (u_n + z_n) + a_j^t, \quad \forall j \in J^H \label{mod:ht:mass} \\
   & 0 \leq v_j \leq V_j, \ \ 0 \leq u_j \leq U_j, \ \ 0 \leq z_j, \quad \forall j \in J^H \label{mod:ht:hyd} \\
   & 0 \leq g_j \leq G_j, \quad j \in J^G\label{mod:ht:ter}\\
   & 0 \leq r_j \leq \tilde{R}_j(\omega^{t}), \quad \forall j \in J^R \label{mod:ht:ren} \\
   & a_j^t = \sum_{l\in J^L} \phi_{j,l} a_j^{t-l} + \tilde{\varepsilon}_j(\omega^{t}), \quad \forall j \in J^H\label{mod:ht:end}
\end{align}

In the above model, the main sources of uncertainty are variable renewable energy (VRE) generation and inflows. The first has a very weak time dependency structure and is considered stagewise independent. The latter, inflows, have a stage-wise dependency that can be modeled using the autoregressive formulation above. The stochastic model of inflows follows the standard Periodic Auto-Regressive (PAR) used by the Brazilian ISO, but here for individualized hydros (in opposition to the aggregated ones used by the ISO) \cite{jardim2001stochastic}. More details of the generation of the joint VRE and inflow scenarios can be obtained in \cite{dias2020modeling}.

The model fits the general form \eqref{mod:mso:beg}--\eqref{mod:mso:end} and can be solved by SDDP.
We describe the result of an SDDP simulation as follows.
Once optimized, we can perform a \textit{simulaton} (see Section \ref{sec:sddp}) and store results as sets vector. For example, $\mathcal{V} = \{\textbf{v}_{t,s}\}_{t\in T,s\in S}$, $\mathit{\Pi} = \{\pi_{t,s}\}_{t\in T,s\in S}$ represent, respectively, volumes and spot prices, (dual variables of the load balance, \eqref{mod:ht:load}), meaning that we have solutions for each sampled scenario $s \in S$ and stage $t \in T$.

The main simplification in this section is to consider the hydro generation function, $\rho$, in \eqref{mod:ht:mass}, as a purely linear function of the turbine flow \cite{rosemberg2021assessing, lohndorf2019modeling, debia2021strategic}, nonlinear models that also consider the dependency of the net-head can also be used, but would require approximations to satisfy the requirements of SDDP \cite{fredo2019assessing}.

\vspace{-3mm}
\section{Strategic agents} \label{sec:sa}

The optimization of an independent agent, or owner, can also be modeled as a multi-stage stochastic optimization problem. In this section, we describe the dynamic model that will considered, introducing features such as the Markovian approach, convexification, and contracts step-by-step. 
\vspace{-3mm}
\subsection{The dynamic model} \label{sec:sa:mso}

The first version of the single price-maker hydro owner agent was modeled and solved with SDDP in \cite{flach2010long}, but it did not consider uncertain bids from other agents. Hence, we present an extension of the technique that handles the uncertainty from other agents' bids with a Markov model. Model \eqref{mod:sa:init}--\eqref{mod:sa:end}, representing the optimization of an individual agent, $i$, is very similar to \eqref{mod:ht:init}--\eqref{mod:ht:end}, the centralized operation model. \eqref{mod:sa:init} represents a Bellman recursion analogous to the one of the centralized dispatch. The main difference in the objective function is that there is a new term to describe the agent's revenue in the given stage $t$ for a random event $\omega^{t}$ as a function of the electricity $e$ produced by the agent. \eqref{mod:sa:load} states that $e$ equals the total generation among all resources of a given agent $i$. Finally, \eqref{mod:sa:mass}--\eqref{mod:sa:end} are almost the same as \eqref{mod:ht:mass}--\eqref{mod:ht:end} but only accounting for generators owned by agent $i$.
\hspace{-1cm}\begin{align}
   &\hspace{-1.7mm}  \tilde{B}^{t}\big(\textbf{v}^{t}_i, \textbf{a}^{[t-1]}_i, \omega^{t}\big) = \min \ \ - \tilde{\Lambda}(e, \omega^{t}) + \sum_{j \in J^G_i} C_j g_j + \notag \\
   &\hspace{-1.7mm}  \hspace{1.3in} \mathbb{E}_{\omega^{t+1}}\left[\tilde{B}^{t+1}\big(\textbf{v}^{t+1}_i, \textbf{a}^{[t]}_i, \omega^{t+1}\big) \right] \label{mod:sa:init}\\
   &\hspace{-1.7mm} s.t.   \notag  \\
   &\hspace{-1.7mm}  e = \sum_{j\in J^G_i} g_j + \sum_{j\in J^H_i} \rho_j(u_j) + \sum_{j\in J^R_i} r_j \label{mod:sa:load}\\
   &\hspace{-1.7mm} v^{t+1}_j = v^t_j - u_j - z_j + \sum_{y\in J^U(j)} (u_y + z_y) + a_j^t,   \forall j \in J^H_i \label{mod:sa:mass}\\
   &\hspace{-1.7mm}  0 \leq v_j \leq V_j, \ \ 0 \leq u_j \leq U_j, \ \ 0 \leq z_j, \forall j \in J^H_i \label{mod:sa:hid}\\
   &\hspace{-1.7mm}  0 \leq g_j \leq G_j, \quad \forall j \in J^H_i \label{mod:sa:ter}\\
   &\hspace{-1.7mm}  0 \leq r_j \leq \tilde{R}_j(\omega^{t}), \quad \forall j \in J^R_i\label{mod:sa:ren}\\
   &\hspace{-1.7mm}  a_j^t = \sum_{l\in J^L} \phi_{j,l} a_j^{t-l} + \tilde{\varepsilon}_j^t(\omega^{t}), \quad \forall j \in J^H_i\label{mod:sa:end}
\end{align}
Although renewable energy generation is still considered stagewise independent and the inflow dependency is still handled by the auto-regressive process of Section \ref{sec:ts-sddp}, the revenue function uncertainty is not restricted to the RHS. Hence,
as highlighted in Section \ref{sec:m-sddp}, the solution of the above model will require the Markovian SDDP, in fact, a mixed auto-regressive (for inflows) and a Markovian (for the revenue) SDDP. The uncertainty representation will be detailed in the next paragraphs.

A second fundamental challenge in this procedure is requiring all hydro plants in the same cascade to belong to a single agent. This hypothesis is also assumed in previous works like \cite{steeger2017dynamic}. Handling different owners in the same river system can be done by considering a water wholesale market besides the electricity market as done in \cite{lino2003bid}. Alternatively, other market designs can be used. Two possible alternative designs are: 1) \textit{hydro slicing}, in which agents own fractions of the cascade and, consequently, can optimize their strategies as if they did not share the cascade, and 2) \textit{virtual hydro reservoirs}, in which all hydro plants of a system are aggregated in a single reservoir which is then split proportionally among all agents. The interested reader is referred to \cite{barroso2012overview} for more details.

\vspace{-3mm}
\subsection{Endogenous spot prices and convexification of the revenue function}\label{sec:sa:rev}
In the case of price-taker agents, their operation can be optimized by considering the revenue function, as in \cite{gjelsvik1999algorithm}:
\begin{align}
    & \tilde{\Lambda}(e, \omega) = \pi(\omega) e \label{mod:sa:revpt}
\end{align}
where $\pi(\omega)$ is a scenario of spot prices that might be obtained from any model. In particular, we could consider spot price scenarios $\mathit{\Pi}$ produced by the \textit{simulation} from Section \ref{sec:ht}.

On the other hand, price-maker agents have the ability to affect spot prices with their energy offers, as in the Stackelberg Equilibrium of Section \ref{sec:st}. In other words, $\pi$ is no longer an exogenous input data, it is a function of the energy offer of agent $i$, i.e., $e$. In this case, price and quantity bids from other agents are considered random variables depending on the uncertainty $\omega$, and are given $\{(P_j(\omega), Q_j(\omega)), j \in I_{-i}\}$. With those bids, we follow the same rationale of \cite{flach2010long, steeger2017dynamic}. Therefore, we write the following model that expresses a simple bid-based dispatch problem:
\begin{align}
  \pi(e, \omega) \in \arg \min_{q,\pi} & \sum_{j \in I_{-i}} P_j(\omega) q_j \label{mod:pi:init}\\
   s.t.&  \sum_{j \in I} q_j = {D} \ \ : \ \pi \label{mod:pi:bal}\\
   & 0 \leq q_i \leq e  \label{mod:pi:limit} \\
   & 0 \leq q_j \leq Q_j(\omega), \quad \forall j \in I_{-i} \label{mod:pi:end}
\end{align}
This is analogous to a market clearing problem, in which a system operator selects the optimal amounts of energy to attain the demand in \eqref{mod:pi:bal}, under the limits \eqref{mod:pi:limit}--\eqref{mod:pi:end} given by the offers of the current agent and the other players. The main result is the spot prices (the dual variable of \eqref{mod:pi:bal}). Now we define:
\begin{align}
    & \tilde{\Lambda}(e, \omega) = \pi(e, \omega) e\label{mod:sa:revpm}
\end{align}

However, this strategic agent revenue function is a non-convex piece-wise linear discontinuous function as detailed in \cite{flach2010long, steeger2017dynamic}. Thus, to satisfy the SDDP convexity requirements, we follow the method proposed in \cite{flach2010long} and represent the convex hull of $\tilde{\Lambda}(e, \omega)$ with respect to $e$ for a fixed $\omega$:
\begin{align}
    \hspace{-0.2cm}chull(\tilde{\Lambda}(e, \omega)) = \max_{\lambda_j \geq 0} & \sum_{j \in J^V} \lambda_j E^R_j(\omega)\label{mod:ch:init}\\
    s.t. &  \sum_{j \in J^V} \lambda_j E^Q_j(\omega) = e,  \sum_{j \in J^V} \lambda_j = 1 \label{mod:ch:end}
\end{align}
where each pair $(E^Q_j(\omega),E^R_j(\omega))$ represents a vertex (in the set of vertices $J^V$) of the convex hull of the hypo-graph of ${\tilde{\Lambda}}(e, \omega)$.
We keep the vertices indexed by $\omega$ to make it clear that they depend on the random variables of the problem. Theoretically, it would be possible to exactly represent the non-convex revenue function with methods like SDDiP \cite{zou2019sddip}. However, this would come with the cost of a considerably larger computational burden. Because we are accounting for many computationally intensive features in this work, we adopted the convex hull representation.

\vspace{-3mm}
\subsection{Accounting for forward contracts in the dynamic model}
As one of the key mechanisms to mitigate market power, we must also be able to represent forward contracts \cite{kelman2001market}. This was not described in the previous MSO models solved by SDDP. However, it is simple to modify the revenue function to consider two additional terms as follows:
\begin{align}
   & \tilde{\Lambda}(e, \omega) = P^F Q^F - \pi(e, \omega) Q^F + \pi(e, \omega) e 
   \label{eq:ctr}
\end{align}
the first term is the fixed revenue of the forward contract, the second represents the energy that must be delivered due to the contract, and the third is the previously represented earnings from the spot market. The constants $P^F$ and $Q^F$ are input data. Consequently, the first term is constant and does not affect the optimization of a risk-neutral (strategic or not) agent. On the other hand, agents are discouraged from generating below contracted quantities (where $e - Q^F < 0$) as they will have to buy energy to match their obligations, thereby losing money. This has an effect that opposes the willingness of price-maker agents to withhold energy (reducing $e$) to increase spot prices.

% they are not decision variables in the optimization problems.

\begin{figure}[!htbp]
\centering
\includegraphics[width=3.4in]{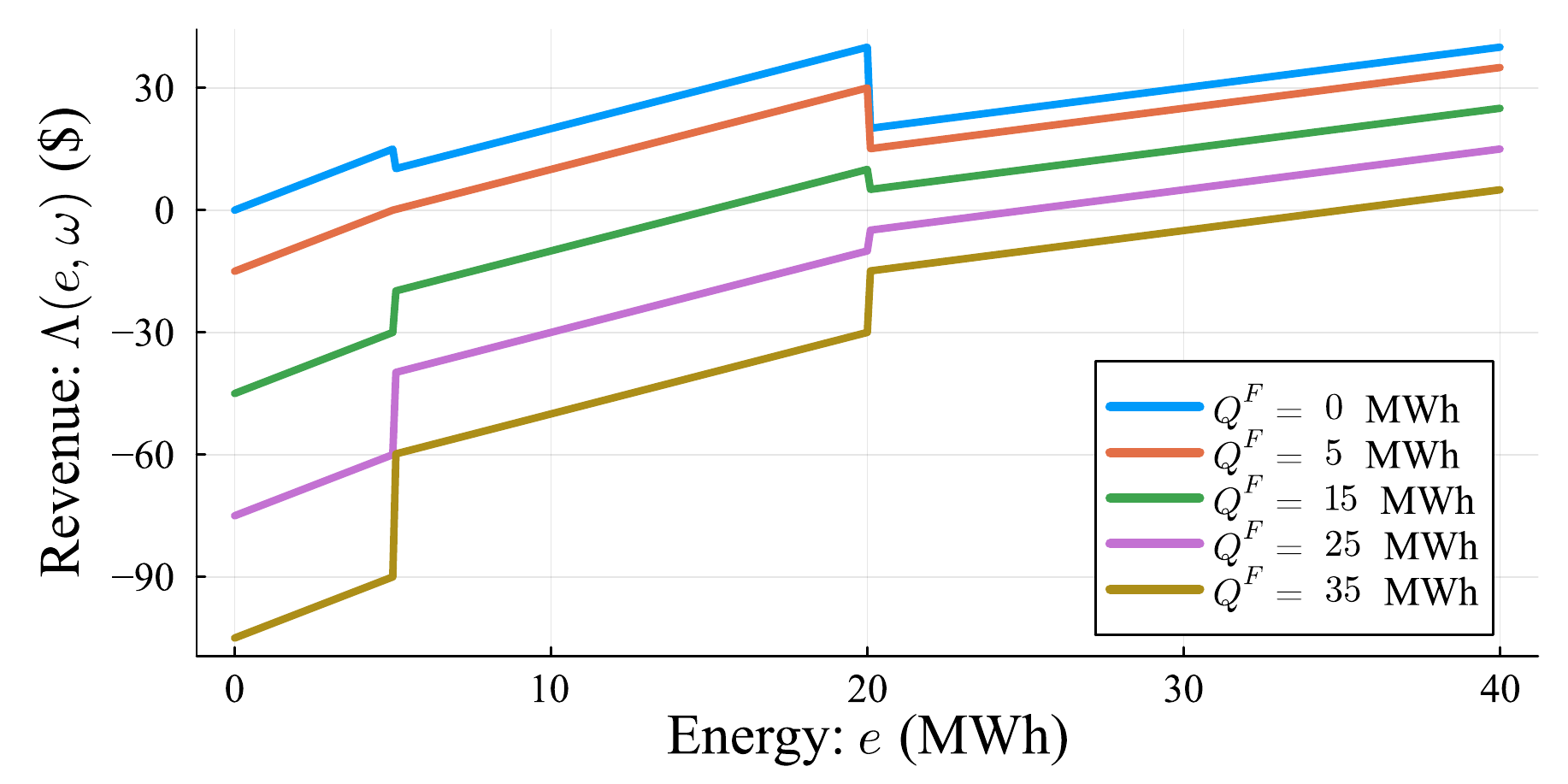}
\caption{Revenue curves, $\tilde{\Lambda}(e, \omega)$, for various values of $Q^F$, and $P^F = 0$.}
\label{fig:sa:seesaw}
\end{figure}

Figure \ref{fig:sa:seesaw} shows revenue curves for a simple case where $P^F=0$\$/MWh, $D = 40$MWh, and we consider $3$ price-quantity offers from the other agents: (3\$/MWh,10MWh), (2\$/MWh,15MWh), (1\$/MWh,20MWh). Note that the total quantity from the sum of the other agents' offers is $45$MWh, which is higher than the demand, and, hence, no deficit occurs even if $e = 0$MWh. This function is also non-convex, but we can represent its convex hull in the optimization problem in the exact same way as the previously described case with no contracts. For instance, we consider $Q^F=0$MWh for which we only have the last term of \eqref{eq:ctr}. In this case, the revenue is a linear function of $e$, with a slope equal to the spot price, $\pi$. The most expensive accepted bid sets the spot price, then a bid $e<5$MWh for the strategic hydro leads to $\pi = 3$\$/MWh. When the strategic hydro offers more than $5$MWh, it displaces the bid previously defining the spot price, and the price drops to the variable cost immediately below, i.e., $2$\$/MWh. So, from this point on, the curve has a smaller slope, and the first break appears. The other curves and breaks follow the same rationale but for different values of $Q^F$ and $e$. For more details, we refer the reader to  \cite{flach2010long}, which only considers the last term in \eqref{eq:ctr}. We emphasize that for different values of $Q^F$, the first term in \eqref{eq:ctr} does not affect this example; as for simplicity, we considered $P^F=0$\$/MWh. However, the second term subtracts a piece-wise constant function from the revenue function as $Q^F$ is constant for each curve and $\pi(e, \omega)$ is a piece-wise constant function. Because $P^F = 0$\$/MWh, the blue curve, which considers no contract, is the highest. In more realistic cases, these curves can cross each other.

\vspace{-4mm}
\subsection{Markovian uncertainty representation}\label{sec:mk}

We remark that in both \cite{flach2010long} and \cite{steeger2017dynamic}, the revenue function is not stochastic since it only considers thermal bids in the form of thermal installed capacities and operating costs and does not require the Markovian SDDP from Section \ref{sec:m-sddp}.

In contrast, we allow stochastic and time-dependent uncertainty in the two cases: 1) spot prices in the price-taker case and 2) energy bids from other agents in the price-maker case.
In both cases, we assume we have scenarios of either spot prices, $\mathit{\Pi}$, or price and quantity bids, $(\mathcal{P},\mathcal{Q})$. Such sets of scenarios are jointly distributed with the inflow, $\mathcal{A}$, and renewable energy, $\mathcal{R}$ scenarios if they were generated by SDDP \textit{simulations}.

Given all these jointly distributed scenarios, we apply the same Expectation Maximization (EM) algorithm as \cite{silva2021data} to obtain the transition probabilities of the Markov model for all stages, $\mathcal{M}$, and the classification of each sample, $s$, of each stage, $t$ in a Markov state.

\vspace{-2mm}
\section{Multiple agents equilibrium simulation} \label{sec:ma}

\begin{figure}[!htbp]
\centering
\includegraphics[width=2.8in, trim={0.cm 0.3cm 0.6cm 0.cm}, clip]{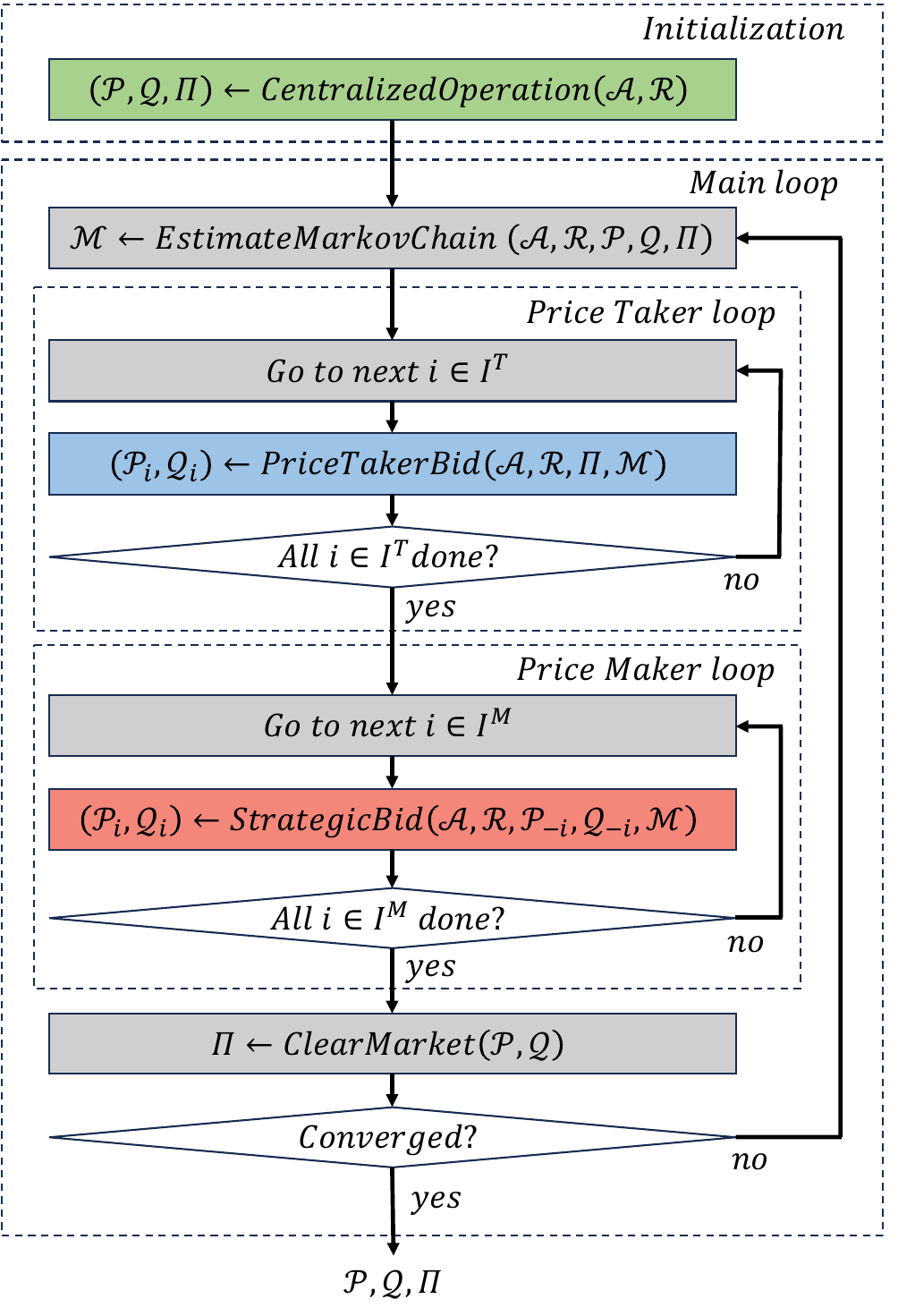}
\caption{Algorithm for multiple agents equilibrium flowchart.}
\label{fig:algo}
\end{figure}

In this section, we combine the methodologies presented in Section \ref{sec:fw} with the ones developed in Sections \ref{sec:ht} and \ref{sec:sa} to describe a novel algorithm to simulate the interaction among multiple price-maker and price-taker agents in a power system with high storage capacity.
The interactions among agents will be in terms of price and quantity bid offers as in a realistic competitive electricity market.
In terms of a game theoretical model,
we have two types of agents, price takers and price makers, whose actions are price and quantity bids, as in Section \ref{sec:nash}.
The complete algorithm is based on the diagonalization procedure of Section \ref{sec:diag} and is depicted in the flowchart of Figure \ref{fig:algo}. All details of the algorithm are presented next.

\vspace{-3mm}
\subsection{Initialization via Centralized Operation}

To start the diagonalization-based algorithm, we need an initial set of actions for all agents. Such actions will be obtained from a \textit{simulation} of the centralized audited-cost-based market design of the power system through the typical cost-minimization model, \eqref{mod:ht:init}--\eqref{mod:ht:end}, from Section \ref{sec:ht} solved with SDDP from Sections \ref{sec:sddp} and \ref{sec:ts-sddp}. Besides the physical description of the system, the main input data for this model are the inflows and renewable generation scenarios, $(\mathcal{A},\mathcal{R})$.
The main results from the simulation procedure are, at a given stage $t$ and scenario $s$:
1) the generation decisions of all units of a generation agent $i$, which lead quantity part of the bid, $Q_{i,t,s}$, and 2) the spot prices of the system ($\pi_{s, t}$) that will be to the price parcel of the bids.
We denote the set of price and quantity bids of an agent $i$ as $(\mathcal{P}_i,\mathcal{Q}_i) = \{P_{i,t,s}, Q_{i,t,s}\}_{t\in T,s\in S}$ and the spot prices as $\mathit{\Pi}$.
We will refer to this procedure as $CentralizedOperation(\mathcal{A},\mathcal{R})$.

\vspace{-3mm}
\subsection{The main loop}

In this step, the actions of each agent will be updated, while the actions of the other ones will be fixed until the convergence criterion is satisfied. The main goal is to simulate the power market in an agent-based fashion. If convergence is strictly attained, we might have reached a Nash equilibrium. However, such an equilibrium might not even exist in this setting. The convergence criterion considered in this work is assuring that the maximum absolute variation of the price and quantity bids of all agents between two subsequent iterations vary less than 1\%  of their values in the centralized operation.

\subsubsection{Markov Chain Estimation} This is the very first part of the main loop because the optimization algorithms presented in the following steps will need a Markov Chain representation of uncertainty of the joint process that includes inflows, renewable generation, bids, and spot prices. The methodology for estimating such Markov Chains, $\mathcal{M}$, was presented in Section \ref{sec:mk}.
We will refer to this procedure as $EstimateMarkovChain(\mathcal{A},\mathcal{R},\mathit{\Pi}, \mathcal{P},\mathcal{Q})$.

\subsubsection{Price-Taker Loop} This loop consists of the self-optimization of all price-taker agents, $i\in I^T$. The model of price takers, from Section \ref{sec:sa}, is given by \eqref{mod:sa:init}--\eqref{mod:sa:end} with the revenue function \eqref{mod:sa:revpt} and it is optimized with a Markovian SDDP from Section \ref{sec:m-sddp}. The uncertainty scenarios are given by the inflows, $\mathcal{A}$, renewable generation, $\mathcal{R}$, and the exogenous spot prices, $\mathit{\Pi}$. For the agent $i$, the \textit{simulation} that succeeds the \textit{policy optimization} will return offered energy, $e_{t,s}$, that together with the input spot prices will define the bid $(\mathcal{P}_i,\mathcal{Q}_i)$. We name this procedure $PriceTakerBid(\mathcal{A},\mathcal{R},\mathit{\Pi}, \mathcal{M})$

\subsubsection{Price-Maker Loop} This loop is analogous to the previous one. The model of each agent $i\in I^M$ is still given by \eqref{mod:sa:init}--\eqref{mod:sa:end}, but the revenue function will be the convex reformulation of the prime-maker revenue function \eqref{mod:ch:init}--\eqref{mod:ch:end}. The model is optimized with a Markovian SDDP from Section \ref{sec:m-sddp}. However, the uncertainty contains the bid from other agents $(\mathcal{P}_{-i},\mathcal{Q}_{-i})$. For the agent $i$, the \textit{simulation} that succeeds the \textit{policy optimization} will return offered energy, $e_{t,s}$, that together with the endogenous spot price, computed from \eqref{mod:pi:init}--\eqref{mod:pi:end}, will define the bid $(\mathcal{P}_i,\mathcal{Q}_i)$. We name this procedure $StrategicBid(\mathcal{A},\mathcal{R},\mathcal{P}_i,\mathcal{Q}_i, \mathcal{M})$

\subsubsection{Market Clearing} This is the last step of the main loop and is responsible for computing update spot prices $\mathit{\Pi}$ by clearing the market for each stage and scenario given bids from all agents, $(\mathcal{P},\mathcal{Q})$.
The model used here is given by  \eqref{mod:pi:init}, \eqref{mod:pi:bal} and \eqref{mod:pi:end}, but with $I_{-i}$ replaced by $I$.
We label this procedure as $ClearMarket(\mathcal{P},\mathcal{Q})$.

\vspace{-2mm}
\section{Computational Experiments with The Brazilian Southeastern System} \label{sec:case}

In this section, we provide quantitative simulated results for the Brazilian power system. We begin with a sensitivity analysis based on the Brazilian Southeast subsystem, which accounts for about $55\%$ of the Brazilian hydro resources and about $50\%$ of the system's total installed capacity. We analyze different market concentrations and contracting schemes based on the procedure developed in Section \ref{sec:ma}. 

\vspace{-3mm}
\subsection{Dataset, infrastructure, and system configuration}\label{sec:case:data}
 
This system was constructed from real data from the Brazilian system expansion scenario for $2025$ (based on the Brazilian market operator data \cite{nw_ccee}) and contains $46$ thermal plants and $21$ hydro plants with a simplified system topology. Note that this is already more than the current $13$ aggregated reservoirs considered in the official model.
The overall hydro installed capacity represents $70\%$ of the system's installed capacity.
To illustrate the algorithm's functioning, first, we consider three main hydro generation companies, one considered as a price-taker agent, and the other two as price makers. The $21$ hydro plants are split into the $3$ aforementioned hydro generation companies, $7$ plants for each one. The three different agents do not share cascades as mentioned in Section \ref{sec:sa}. All thermal plants are considered individual price-taker agents. We considered a single load block for the sake of simplicity, that is, a constant demand for the entire month.
The following simulations were carried out in a $5$-year monthly horizon, that is, $60$ stages. We considered $1000$ sampled scenarios.
The large-scale simulations were performed in a \textit{PSRCloud} cluster of $8$ servers, each one with $64$ cores and $128$GB of RAM. Each simulation took approximately $3$ hours. All algorithms were coded in the Julia language \cite{bezanson2017julia},
and the optimization models were coded with JuMP \cite{Lubin2023} and solved with the Xpress solver\cite{xpress}.

\vspace{-3mm}
\subsection{Impact of Market Concentration on Market Power}

\begin{figure*}[!htbp]
   \centering
   \subfloat[]{\includegraphics[width=3.4in, trim={0cm 0.5cm 0cm 0.5cm}, clip]{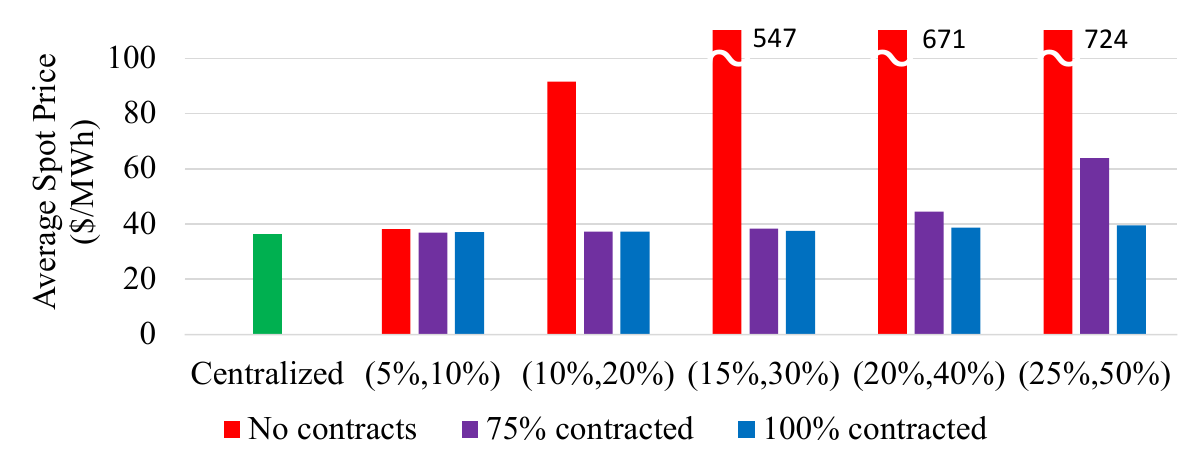}%
   \label{fig:se:spot}}
   \hfil
   \subfloat[]{\includegraphics[width=3.4in, trim={0cm 0.5cm 0cm 0.5cm}, clip]{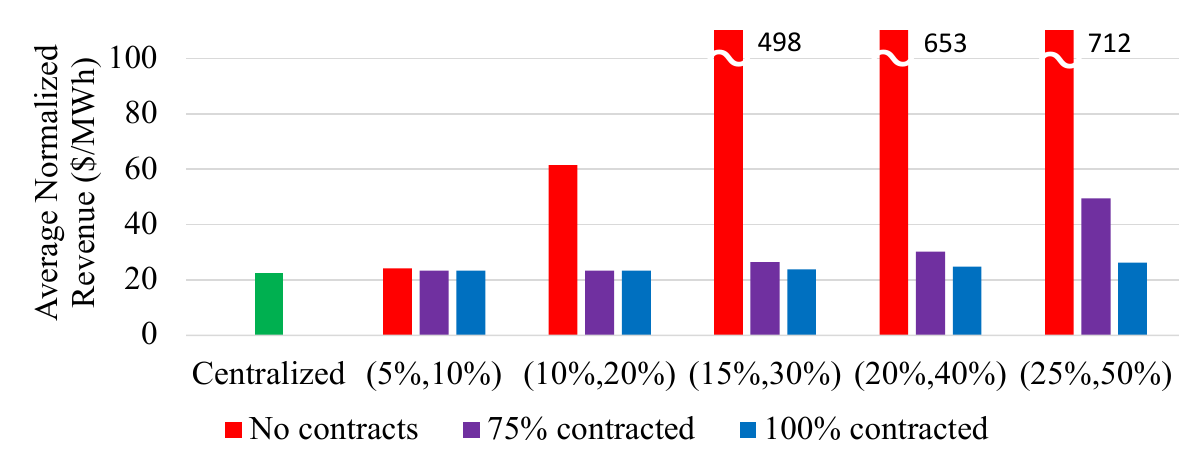}%
   \label{fig:se:rev}}
   % \centering
   \hfil
   \subfloat[]{\includegraphics[width=3.4in, trim={0cm 0.45cm 0cm 0.5cm}, clip]{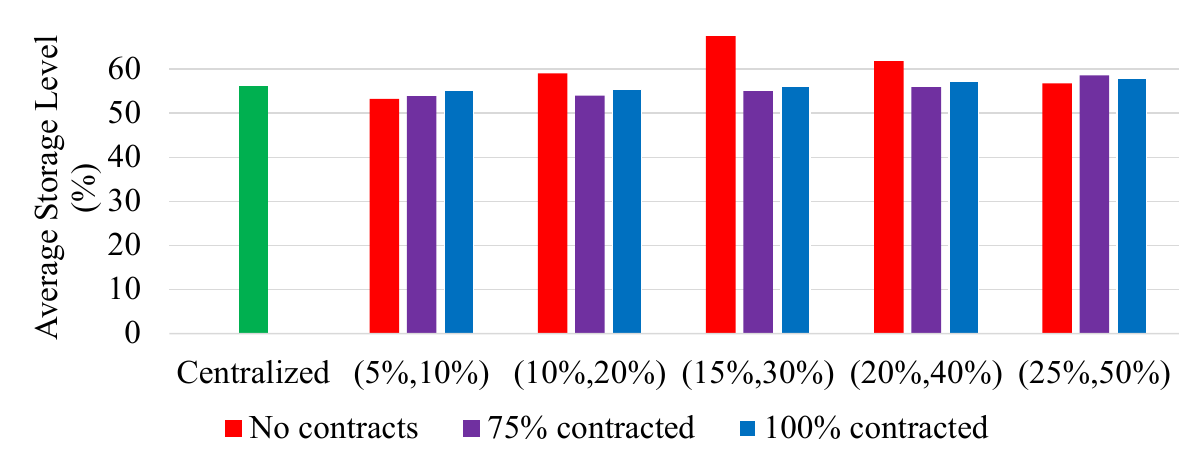}%
   \label{fig:se:stor}}
   \hfil
   \subfloat[]{\includegraphics[width=3.4in, trim={0cm 0.45cm 0cm 0.5cm}, clip]{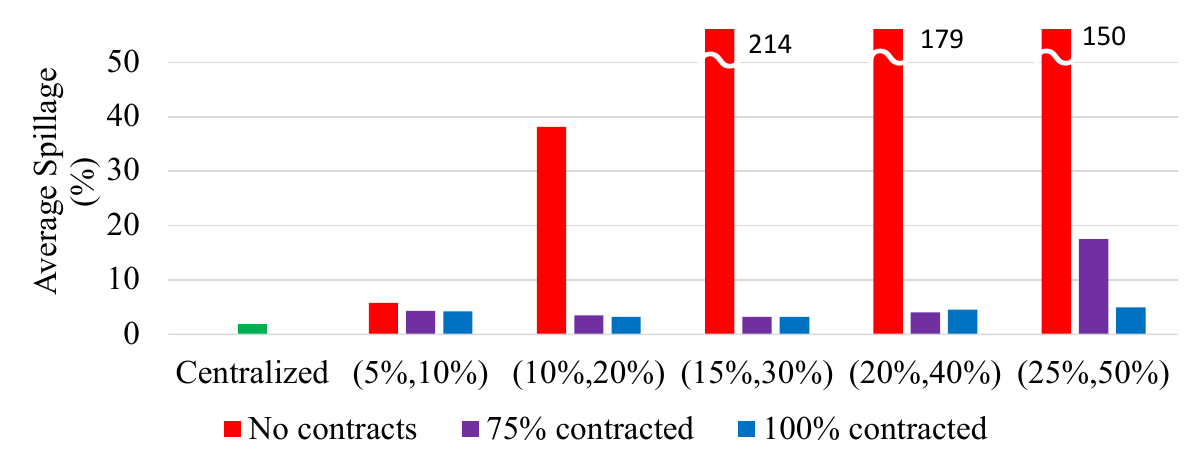}%
   \label{fig:se:spill}}
   % \centering
   \caption{Results for simulations of Brazil's Southeast under different market concentrations.
   All plots have the same template and contain average values of: (a) Spot Prices, (b) Normalized Revenue, (c) Spillage, and (d) Storage Level. Averages are with respect to all stages and scenarios. Spillage and Storage $\%$ are with respect to the maximum amount of water that can be stored in the system.
   All plots include a first bar (in green) with the result for the centralized case.
   In parentheses in the horizontal axis, we have the percent share of each price-maker agent, defining the database used.
   For each database, we have 3 bars, one for a case with no contracts (in red) and another two for cases with 75\% \textit{of contracting level} (in purple) and 100\% \textit{of contracting level} (in blue). \\ Note that some red bars in (a), (b), and (d) are too high and do not fit the plot area, so they contain the values beside the top of the bar.}
   \label{fig:se}
\end{figure*}

In this section, we present five controlled simulation experiments to understand the effect of different market concentrations on market power abuse. We adjust the hydro plants so that we have different databases defined by tuples $({share}_1\%, {share}_2\%)$. In each case, ${share}_1$ and ${share}_2$ represent the shares (in percentage) of the hydro system belonging to each of the two price-maker companies (agents). The number of plants of each hydro agent is fixed as mentioned in Section \ref{sec:case:data}.
We multiplied and installed capacities and reservoir sizes by scalars to adjust the hydro system size of each agent in each database. The remainder of the hydro system is assumed to be of price-taker agents.

\begin{figure}[!htbp]
\centering
\includegraphics[width=3.4in]{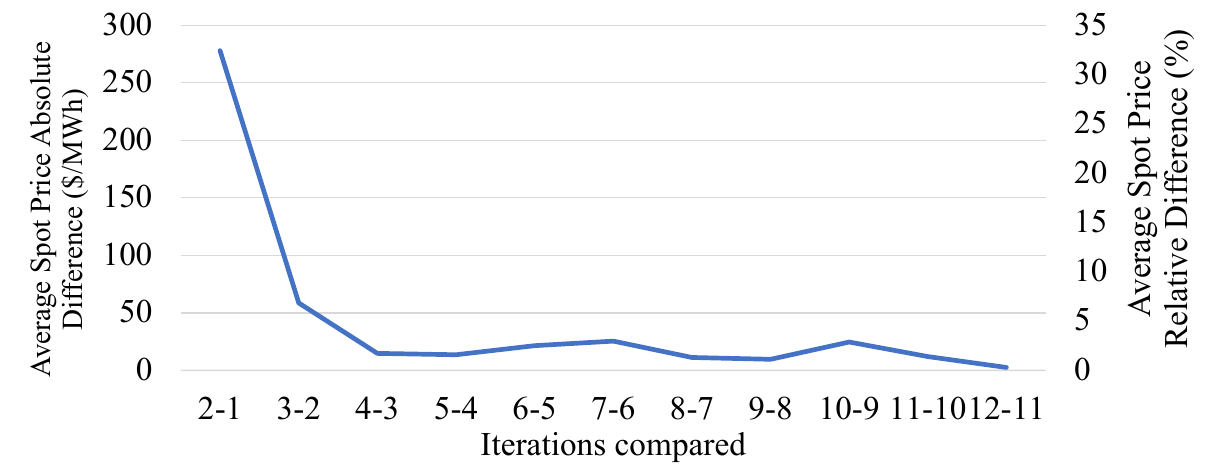}
\caption{Convergence profile from Brazil's southeast with $(25\%,50\%)$ and \textit{no contracts}. See Figure \ref{fig:se}. Average spot price absolute and relative differences between two consecutive iterations: $(i+1)$--$(i)$.}
\label{fig:se:conv:diff}
\end{figure}

Figure \ref{fig:se:conv:diff} presents a typical convergence profile of the proposed iterative method. These profiles were taken from the case study with price-maker shares of $(25\%, 50\%)$ and no contracts. Figure \ref{fig:se:conv:diff} shows both average absolute and relative differences between two consecutive iterations. First, absolute (or relative) differences are computed for each stage and scenario, and then, they are averaged. 

Figure \ref{fig:se} contains results from this first study and from the next study for better understanding. For now, we shall focus on the centralized cases (in green) and the multiple cases with no contracts (all in red). In Figure \ref{fig:se:spot}, we show the average spot prices as a function of the market distributions in a bar plot. We add the first bar with the spot price of the centralized dispatch. As expected, the spot prices rapidly grow as the concentration increases. Figure \ref{fig:se:rev} depicts the captured price, i.e., the average value of energy sold, which is the total revenue of each agent divided by its total generated energy. We present results in such normalized forms to compare the results of different-sized agents. Figures \ref{fig:se:stor} and \ref{fig:se:spill} show additional results with average reservoir levels of the price-maker agents throughout the study period and the average spillage. Interestingly, market power abuse is mostly characterized by excess spillage, whereas, on average, the total reservoir levels do not vary much compared to the benchmark (centralized dispatch with audited costs).

\vspace{-4mm}
\subsection{Impact of Contracts on Market Power}

In this section, we repeat the above analyses but with contracts. To come up with monthly contract quantities, we take average generation values from the centralized dispatch and prices from the average spot prices in the centralized dispatch. This will stimulate the agents to produce energy and reduce the spot prices since being short in the contract together with high spot prices will lead the agent to have expenses due to the second term in \eqref{eq:ctr}.

Now we look again to Figures \ref{fig:se:spot}, \ref{fig:se:rev}, \ref{fig:se:stor} and \ref{fig:se:spill} described in the previous section but we also will focus on the cases of agents that are $75\%$ (purple) and $100\%$ (blue) contracted, respectively. In the case studies, we can clearly see that the contracts almost completely eradicated the market power, and the resulting spot prices are very close to the ones obtained in the centralized dispatch. At the same time, the captured revenues and spillage levels were shrunk to values close to the ones observed on the centralized dispatch with audited costs. By analyzing Figure \ref{fig:se:spill}, we can see that spillage is a relevant marker for market power abuse.

\vspace{-2mm}
\section{Case Study: The Brazilian Power System}\label{sec:br}

In this section, we provide results for a case study considering the complete Brazilian power system. We contrast results from simulations with and without contracts and provide a first quantitative measurement of the market power abuse potential in this system. 

\vspace{-3mm}
\subsection{Dataset, infrastructure, and system configuration}
This dataset was created from the original data of the Brazilian system obtained from the Brazilian market operator webpage \cite{nw_ccee}. The dataset contains
$137$ thermal plants representing $18\%$ of the installed capacity, $364$ renewable energy plants (wind and solar) representing $19\%$ of the installed capacity and $32$ hydropower plants representing the remaining $63\%$ of the installed capacity. Here, we improve the load duration curve accuracy and consider $5$ load blocks \cite{roh2009market, wogrin2014new} instead of just one (from the previous section) to better represent peak demand hours, $5$ is more than the $3$ currently used in the Brazilian official planning tools \cite{nw_ccee}.
Like in previous studies, we consider $5$ years ($60$ monthly stages), $1000$ scenarios, and the same \textit{PSRCloud} cluster configuration.

In this study, we consider an approximation of the real market concentration in the system. We represent three companies as price-maker agents, with respectively $32\%$, $9\%$, and $7\%$ of the resources in terms of firm energy representing the three largest generation owners. The firm energy metric is used as hydro and renewable capacity is not fully available all the time, so this is a proxy of their average available generation, see \cite{faria2009allocation}. Meanwhile, the remaining $55\%$ of the generators are distributed among one price-taker hydro company, and each thermoelectric generator is considered another small price-taker agent. For completeness, we computed the Herfindahl-Hirschman Index (HHI) of the system, which led to the value of $0.125$, which is below $0.15$, a number that typically represents a small concentration. However, as will be empirically demonstrated in the next study, there is clear room for market power in the system in question. We believe that it is yet another indication of the limitations of HHI that were already highlighted by \cite{borenstein1999market, kelman2001market}.

Each simulation from this study took approximately $11$ hours using the same code and computational infrastructure of the previous study.

\vspace{-3mm}
\subsection{Market Power in the Brazilian Power System as a Function of the Contracting Level}

\begin{figure}[!h]
\centering
\includegraphics[width=3.4in, trim={0.5cm 0cm 0.5cm 0cm}, clip]{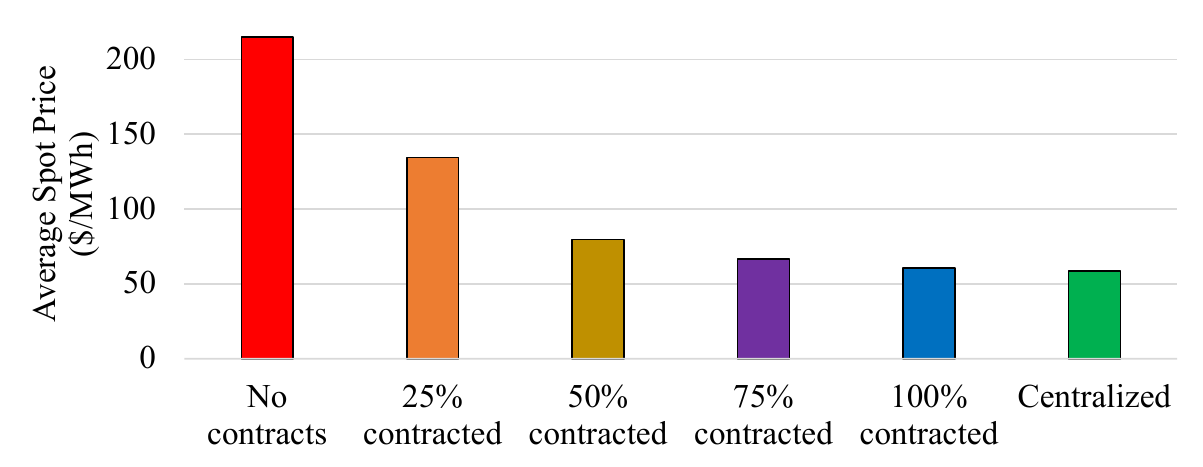}
\caption{Average Spot Price (\$/MWh) for the Brazilian system under different contracting levels. Averages with respect to stages and scenarios.}
\label{fig:br:spot}
\end{figure}

\begin{figure}[!h]
% \centering
\includegraphics[width=3.4in, trim={0.6cm 0cm 0.6cm 0cm}, clip]{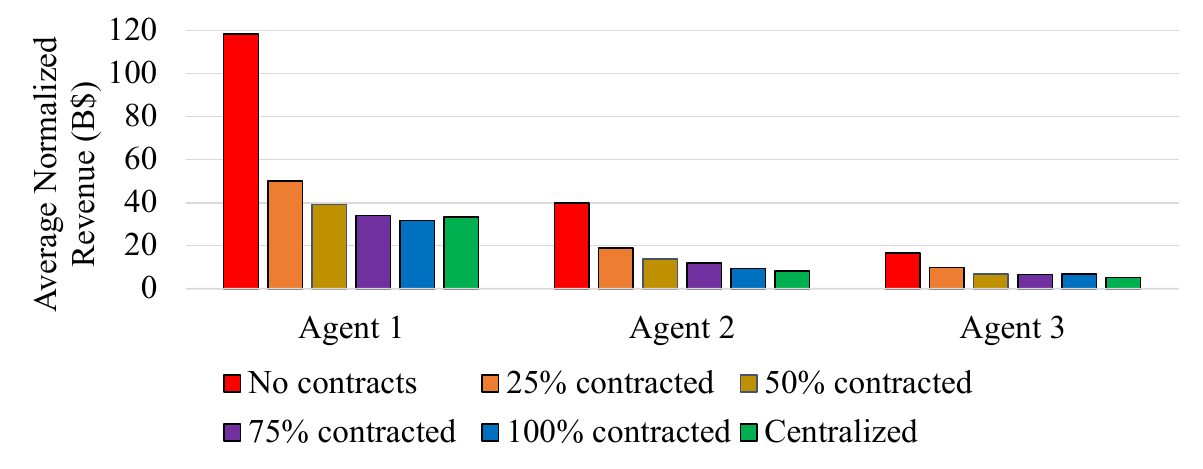}
\caption{Average normalized revenue (\$/MWh) for each of the $3$ price-maker agents under different contracting levels. Averages with respect to stages and scenarios.}
\label{fig:br:rev}
\end{figure}

First, we simulated the centralized operation and the bid-based market without contracts. Then, we considered agents with multiple contracting levels in the bid-based market. In particular, we analysed the cases of: $25\%$, $50\%$, $75\%$ and $100\%$
of contracting level. The key results from each of the $6$ simulations are depicted in Figures \ref{fig:br:spot} and \ref{fig:br:rev}.

Figure \ref{fig:br:spot} shows the spot price under the above $6$ conditions. We can see the progression from the case with the highest spot prices, where no contracts are considered, until the fully contracted case, where the average spot price is much closer to the one obtained from the centralized operation with audited costs. Given the number of resources allocated to the price-maker agents, $45\%$, we note that a $75\%$ contracting level effectively reduces the gap to the centralized operation. This closely follows what we have seen in the previous section with the southeastern system. In such case, $75\%$ was also effective in containing market power in the case where the first price-maker agent had $15\%$ and the second price maker had $30\%$ of the resources, resulting in a $45\%$ share of the system in the hands of price-maker agents. Finally, we note that even low contracting levels, such as $25\%$, significantly reduce the final average spot prices. So, we highlight the empirical evidence in our case study that the typical contracting levels observed in practice, ranging from $75\%$ to $100\%$, can be seen as a relevant instrument to prevent market power abuse in the Brazilian power system. 

Figure \ref{fig:br:rev} presents the overall revenue, including operation and contract costs and spot and contract revenue throughout the $5$-year horizon for each of the three price-maker agents. Overall, the figure reproduces information already verified in Figure \ref{fig:br:spot} with excessive revenues for all agents for the bid-based cases with no contracts. Notably, the $100\%$ contracting case leads to revenues that are very similar to the ones from the centralized dispatch. Interestingly, $75\%$ seems to be a relevant threshold value as it induces behaviors resulting in revenues and spot prices relatively close to the ones obtained in the centralized case. This observation is also relevant because when generators are risk-averse, it is expected that their optimal forward involvement should be slightly below the 100\% case (see \cite{maier2016risk} and references therein).

\begin{figure}[!h]
% \centering
\includegraphics[width=3.4in, trim={0.5cm 0cm 0.6cm 0cm}, clip]{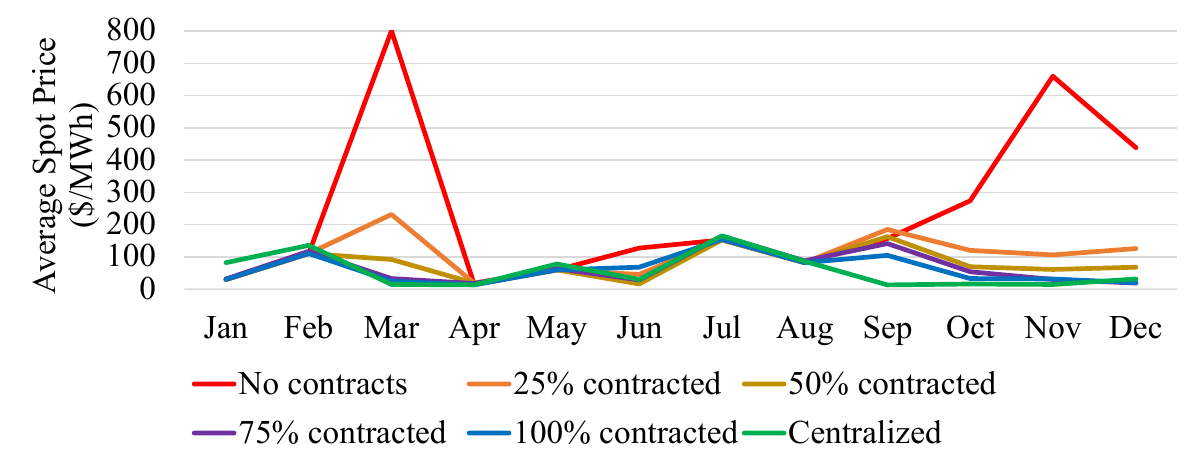}
\caption{Spot price annual profile (MWh).}
\label{fig:br:spot_year}
\end{figure}

Figure \ref{fig:br:spot_year}
presents the seasonal average spot price profile. Spot price changes are driven by both withheld and spilled hydro water. 
Noticeable, the spot-price peaks caused by the market power abuse are concentrated at the end of the dry season, when reservoirs are at their lowest levels (September-December), and in the month with the highest inflow (March), in which strategic agents can easily spill.

\vspace{-2mm}
\section{Conclusion} \label{sec:con}

Based on previous works on long-term hydrothermal power markets, we combined three key methodologies to develop a new and effective market simulator, namely: 1) SDDP applied to the centralized hydrothermal dispatch to initialize and benchmark the process, 2) a multi-stage bilevel stochastic model to model the strategic bidding strategy for eachprice-maker agent, and 3) an iterative diagonalization-based method to simulate the interactions among agents in a market equilibrium. In contrast to the existing literature on the subject, with our new algorithm, we could jointly consider multiple reservoirs (32), multiple stages (60) and scenarios (1,000), and multiple price-maker agents (3). Therefore, a realistic large-scale case study based on the Brazilian power system (one of the largest in the world) could be addressed, and relevant insights could be obtained. 

On the other hand, some simplifications were necessary. In particular, the hydro generation function was considered as a purely linear function of the turbine flow (a standard assumption commented in Section \ref{sec:ht}); agents did not share cascades, as mentioned in Section \ref{sec:sa:mso}; revenue curves were convexified in the SDDP policy optimization, as described in Section \ref{sec:sa:rev}.
We believe that these simplifications can be explored in future works and lead to relevant new insights. So, we believe our simulator and methodology pave the way for relevant new studies.

Within the limitations of the presented computational experiments and case study, which include all assumptions of the proposed model and the specific data, the results and analyses carried out in this work allow us to convey the following concluding remarks for the Brazilian power system:

\begin{enumerate}
    \item We found evidence that the market concentrations observed in the current Brazilian electricity market are capable of enabling market power abuse.
    \item The spillage level is a relevant market for market power abuse.
    \item Contracts are relevant instruments to prevent market power abuse in this country. Our results indicate a significant increase in long-term average spot prices due to market power abuse in cases of lower contracting levels. 
\end{enumerate}

We highlight that the developed market analysis tool can help regulators and market/system operators in market monitoring activities, expansion studies, and understanding the impact of new market rules and instruments. For instance, based on the above remarks, low contracting levels combined with high spillages can be relevant markers of potential market power abuse. On the algorithmic side, the method is suitable for parallel computing and can be extended to consider other physical and regulatory details. Finally, we showcase that a real-world problem formulated as a multi-stage bilevel stochastic problem can be solved under reasonable assumptions and provide relevant insights to decision-makers and regulators. We highlight the following as relevant research avenues: 1) studying the effects and incentives of agents' risk aversion on the equilibrium, 2) studying generators' incentives to contract through a joint spot and long-term forward market equilibrium, and 3) analyzing the impact of different forecast model accuracies on the equilibrium, representing the diversity of views among agents.

\vspace{-4mm}

\bibliographystyle{ieeetr}
\bibliography{ref}

\vfill

\end{document}